\documentclass[10pt]{article}
\usepackage{color}
\usepackage{amssymb}
\usepackage{amsthm,array,amssymb,amscd,amsfonts,latexsym, url}
\usepackage{amsmath}

\newtheorem{theo}{Theorem}[section]
\newtheorem{prop}[theo]{Proposition}
\newtheorem{lemm}[theo]{Lemma}

\newtheorem{rema}[theo]{Remark}

\newtheorem{variant}[theo]{Variant}
\newtheorem{conj}[theo]{Conjecture}

\title{Remarks on curve classes on rationally \\ connected varieties}
\author{Claire Voisin
\\CNRS, Institut de math\'{e}matiques de Jussieu}
\date{}

\newcommand{\cqfd}
{%
\mbox{}%
\nolinebreak%
\hfill%
\rule{2mm}{2mm}%
\medbreak%
\par%
}
\newfont{\gothic}{eufb10}
\begin{document}
\maketitle

\begin{flushright}
{\it To Joe Harris, on his 60th birthday}
\end{flushright}

\section{Introduction}
Let $X$ be a smooth complex projective variety. Define
\begin{eqnarray} Z^{2i}(X)=\frac{{\rm Hdg}^{2i}(X,\mathbb{Z})}{H^{2i}(X,\mathbb{Z})_{alg}},
\end{eqnarray}
where ${\rm Hdg}^{2i}(X,\mathbb{Z})$ is the space of integral Hodge classes on $X$ and
$H^{2i}(X,\mathbb{Z})_{alg}$ is the subgroup of $H^{2i}(X,\mathbb{Z})$ generated by classes of codimension $i$ closed algebraic subsets of $X$.

These groups measure the defect of the Hodge conjecture for integral Hodge classes, hence they are trivial for $i=0$, $1$ and
$n={\rm dim}\,X$, but in general
 they can be nonzero by \cite{AH}. Furthermore they are torsion if the Hodge conjecture for {\it rational} Hodge classes on $X$ of degree $2i$ holds.
In addition to the previously mentioned case, this happens when $i=n-1,\,n={\rm dim}\,X$, due to the Lefschetz theorem on $(1,1)$-classes and the hard Lefschetz isomorphism (cf. \cite{voisinjapjmath}).
We will call classes in ${\rm Hdg}^{2n-2}(X,\mathbb{Z})$ ``curve classes'',
as they are also degree $2$ homology classes.

Note that the Koll\'ar counterexamples (cf. \cite{kollar}) to the integral Hodge conjecture
already exist for  curve classes (that is degree $4$ cohomology classes in this case)
on projective threefolds, unlike the Atiyah-Hirzebruch examples which work for degree $4$
integral Hodge  classes in higher dimension.

It is remarked in \cite{soulevoisin}, \cite{voisinjapjmath} that
the two groups
$$Z^{4}(X),\,\,Z^{2n-2}(X),\,n:={\rm dim}\,X$$
are birational invariants. (For threefolds, this is the same group, but not in higher
dimension.) The nontriviality of these birational invariants for rationally connected
varieties is asked in \cite{voisinjapjmath}. Still more interesting is the nontriviality
of these invariants for unirational varieties, having in mind the
L\"{u}roth problem (cf. \cite{beauville}, \cite{B}, \cite{CG}).

Concerning the group $Z^4(X)$,  Colliot-Th\'{e}l\`{e}ne and the
author proved in \cite{ctv}, building on the work of
Colliot-Th\'{e}l\`{e}ne and Ojanguren \cite{colliotojanguren}, that
it can be nonzero for unirational varieties starting from dimension
$6$. What happens in dimensions $5$ and $4$ is unknown (the four
dimensional case being particularly challenging in our mind), but in
dimension $3$, there is the following result  proved in
\cite{voisinuniruled}:
\begin{theo}\label{theouniruled} (Voisin 2006) Let $X$ be a smooth projective threefold
 which is either uniruled or Calabi-Yau. Then the group $Z^4(X)$ is equal to $0$.
\end{theo}

This result, and in particular the Calabi-Yau case, implies that the
group $Z^{6}(X)$ is also $0$ for a Fano fourfold $X$ which admits a
smooth anticanonical divisor. Indeed, a smooth anticanonical divisor
$j: Y\hookrightarrow X$ is a Calabi-Yau threefold, so that we have
$Z^4(Y)=0$ by Theorem \ref{theouniruled} above. As
$H^2(Y,\mathcal{O}_Y)$, every class in $H^4(Y,\mathbb{Z})$ is a
Hodge class, and it follows that
$H^4(Y,\mathbb{Z})=H^4(Y,\mathbb{Z})_{alg}$. As the Gysin map
$j_*:H^4(Y,\mathbb{Z})\rightarrow H^6(X,\mathbb{Z})$ is surjective
by the Lefschetz theorem on hyperplane sections, it follows that
$H^6(X,\mathbb{Z})=H^6(X,\mathbb{Z})_{alg}$, and thus $Z^6(X)=0$.

In the paper \cite{hoeringvoisin}, it was proved more generally that if
$X$ is any Fano fourfold, the group $Z^6(X)$ is trivial. Similarly, if $X$ is a Fano fivefold on index $2$,
the group $Z^8(X)$ is trivial.

These results have been generalized to higher dimensional Fano
manifolds of  index $n-3$ and dimension $\geq8$  by  Enrica Floris \cite{floris} who
proves the following result:
\begin{theo} Let $X$ be a Fano manifold over $\mathbb{C}$ of dimension $n\geq 8$ and index $n-3$. then the
group $Z^{2n-2}(X)$ is equal to $0$: Equivalently, any integral cohomology class of degree $2n-2$ on $X$ is algebraic.
\end{theo}

The purpose of this note is to provide a number of evidences for the vanishing of
the group $Z^{2n-2}(X)$, for any rationally connected variety over $\mathbb{C}$. Note that in this case, since
$H^2(X,\mathcal{O}_X)=0$, the Hodge structure on $H^2(X,\mathbb{Q})$ is trivial, and so is the
Hodge structure on $H^{2n-2}(X,\mathbb{Q})$, so that $Z^{2n-2}(X)=H^{2n-2}(X,\mathbb{Z})/H^{2n-2}(X,\mathbb{Z})_{alg}$.
We will first prove the following two results.
\begin{prop} \label{propdefoinv} The group $Z^{2n-2}(X)$ is locally deformation invariant
for rationally connected manifolds $X$.
\end{prop}
Let us  explain the meaning of the statement. Consider a smooth projective morphism $\pi:\mathcal{X}\rightarrow B$
between connected quasi-projective complex varieties, with $n$ dimensional
fibers. Recall from \cite{komimo} that if one fiber
$X_b:=\pi^{-1}(b)$ is rationally connected, so is every fiber. Let us  endow everything with the usual topology.
Then the sheaf $R^{2n-2}\pi_*\mathbb{Z}$ is locally constant on $B$. On any Euclidean open set
$U\subset B$ where this local system  is trivial, the group
$Z^{2n-2}(X_b),\,b\in U$ is the finite quotient of the {\it constant} group $H^{2n-2}(X_b,\mathbb{Z})$
by its subgroup $H^{2n-2}(X_b,\mathbb{Z})_{alg}$.
To say that $Z^{2n-2}(X_b)$ is locally constant means that on open sets $U$ as above, the subgroup
$H^{2n-2}(X_b,\mathbb{Z})_{alg}$ of the constant group  $H^{2n-2}(X_b,\mathbb{Z})$ does not depend on $b$.

It follows from the above result that the vanishing of the group $Z^{2n-2}(X)$ for
$X$ a rationally connected  manifold reduces to the similar statement for $X$  defined over a number field.

Let us now define an $l$-adic analogue $Z^{2n-2}(X)_l$ of the group
$Z^{2n-2}(X)$ (cf. \cite{ctk}, \cite{ctz}). Let $X$ be a smooth
projective variety defined over a field $K$ which in the sequel will
be either a finite field or a number field. Let $\overline{K}$ be an
algebraic closure of $K$. Any cycle $Z\in CH^s(X_{\overline{K}})$ is
defined over a finite extension of $K$. Let  $l$ be a prime  integer
different from  $p={\rm char}\,K$ if $K$ is finite. It follows that
the cycle class
$$cl(Z)\in H^{2s}_{et}(X_{\overline{K}},\mathbb{Q}_l(s))$$
is invariant under a finite index subgroup of ${\rm Gal}(\overline{K}/K)$.

 Classes satisfying this property are called Tate classes. The Tate conjecture for
finite fields asserts the following:
\begin{conj}\label{tateconj} (cf. \cite{milne} for a recent account) Let
$X$ be smooth and projective over a finite field $K$.  The cycle class
map gives for any $s$ a surjection
$$cl: CH^{2s}(X_{\overline{K}})\otimes\mathbb{Q}_l\rightarrow H^{2s}(X_{\overline{K}},\mathbb{Q}_l(s))_{Tate}.$$
\end{conj}

Note that the cycle class defined on $CH^s(X_{\overline{K}})$ takes in fact values
in $H^{2s}(X_{\overline{K}},\mathbb{Z}_l(s))$, and more precisely in the subgroup
$H^{2s}(X_{\overline{K}},\mathbb{Z}_l(s))_{Tate}$ of classes invariant under
a finite index subgroup of ${\rm Gal}(\overline{K}/K)$. We thus get for each $i$
a morphism
$$cl^i: CH^{2i}(X_{\overline{K}})\otimes\mathbb{Z}_l\rightarrow H^{2i}(X_{\overline{K}},\mathbb{Z}_l(i))_{Tate}.$$

We can thus introduce the following variant of the groups $Z^{2i}(X)$:
$$Z^{2i}_{et}(X)_l:=H^{2i}_{et}(X_{\overline{K}},\mathbb{Z}_l(i))_{Tate}/{\rm Im}\,cl^i.$$

An argument similar to the one used for the proof of Proposition
\ref{propdefoinv}  will  lead to  the following result:
\begin{prop} \label{spectheo} Let $X$ be a smooth rationally connected variety  defined over a number field
$K$, with ring of integers $\mathcal{O}_K$. Assume given a projective model
$\mathcal{X}$ of $X$ over ${\rm Spec}\,\mathcal{O}_K$. Fix a prime integer $l$. Then
except for finitely many $p\in {\rm Spec}\,\mathcal{O}_K$, the group
$Z^{2n-2}_{et}(X)_l$ is isomorphic to the group
$Z^{2n-2}_{et}(X_p)_l$.
\end{prop}

In the course of the paper, we will also consider  variants $Z^{2n-2}_{rat}(X)$, resp. $Z^{2n-2}_{et,rat}(X)_l$
of the groups $Z^{2n-2}(X)$, resp. $Z^{2n-2}_{et}(X)_l$, obtained by taking the quotient of the group
of integral Hodge classes (resp. integral $l$-adic Tate classes) by the subgroup generated by classes of
{\it rational} curves. This variant is suggested by Koll\'ar's paper (cf. \cite[Question 3, (1)]{kollarportug}).
By the same arguments, these groups are also deformation and specialization invariants for rationally
connected varieties.

Our  last result is conditional but it strongly suggests the
vanishing of the group $Z^{2n-2}(X)$ for $X$ a smooth rationally
connected variety over $\mathbb{C}$. Indeed, we will prove  using
the main result of  \cite{schoen} and the two propositions above the
following consequence of Theorem \ref{spectheo}:
\begin{theo} \label{theoapplischoen} Assume  Tate's conjecture \ref{tateconj} holds for  degree $2$ Tate classes on smooth projective surfaces defined over a finite field. Then  the group
$Z^{2n-2}(X)$ is trivial  for  any   smooth rationally connected variety
 $X$ over $\mathbb{C}$.
\end{theo}
{\bf Thanks.} {\it I thank the organizers of the beautiful conference ``A celebration of algebraic geometry'' for inviting me there.
 I also thank Jean-Louis Colliot-Th\'{e}l\`{e}ne, Olivier Debarre and J\'anos Koll\'ar for useful discussions.

It is a pleasure to dedicate this note to Joe Harris, whose influence on the subject of rational curves on algebraic
varieties  (among other topics!) is invaluable.}

\section{Deformation and specialization invariance}
{\bf Proof of Proposition \ref{propdefoinv}.}
We first observe that, due to the fact that relative Hilbert schemes parameterizing curves in the fibers of $B$
are a countable union of varieties which are projective over $B$, given a simply connected open set $U\subset B$ (in the classical topology of $B$),
and a class $\alpha\in \Gamma(U,R^{2n-2}\pi_*\mathbb{Z})$ such that $\alpha_t$ is algebraic for $t\in V$, where
 $V$ is a smaller nonempty open set $V\subset U$,
then $\alpha_t$ is algebraic for any $t\in U$.

To prove the deformation invariance, we  just need using the  above observation  to prove the following:

\index{integral Hodge classes}
\begin{lemm} \label{le4sep} Let $t\in U\subset B$, and let $C\subset X_t$ be a curve and
let $[C]\in H^{2n-2}(X_t,\mathbb{Z})\cong \Gamma(U,R^{2n-2}\pi_*\mathbb{Z})$ be its cohomology class. Then
the class $[C]_s$ is algebraic for $s$ in a neighborhood of $t$ in
$U$.
\end{lemm}

{\bf Proof of Lemma \ref{le4sep}.} By results of \cite{komimo},
there are rational curves $R_i\subset X_t$ with ample normal bundle
which meet $C$ transversally at distinct points, and with arbitrary
tangent directions at these points. We can choose an arbitrarily
large number $D$ of such curves with generically chosen tangent
directions at the attachment points. We then know by \cite[{\S}2.1]{GHS}
that the curve $C'=C\cup_{i\leq D} R_i$ is smoothable in $X_t$ to
 a smooth unobstructed curve $C''\subset X_t$,
that is $H^1(C'',N_{C''/X_t})=0$. This curve $C''$ then deforms with
$X_t$ (cf. \cite{kodaira}, \cite[II.1]{kollarbook}) in the sense
that the morphism from the deformation of the pair $(C'',X_t)$ to
$B$ is smooth, and in particular open. So there is a neighborhood of
$V$ of $t$ in $U$ such that for  $s\in V$, there is a curve
$C''_s\subset X_s$ which is a deformation of $C''\subset X_t$. The
class $[C''_s]=[C'']_s$ is thus algebraic on $X_s$. On the other
hand, we have
$$[C'']=[C']=[C]+\sum_i[R_i].$$
As the $R_i$'s are rational curves with positive normal bundle, they
are also unobstructed, so that the classes $[R_i]_s$ also are
algebraic on $X_s$ for $s$ in a neighborhood of $t$ in $U$. Thus
$[C]_s=[C'']_s-\sum_i[R_i]_s$ is algebraic on $X_s$ for $s$ in a
neighborhood of $t$ in $U$. The lemma, hence also the proposition,
is proved.

\cqfd
\begin{rema}\label{remarat}{\rm  There is an interesting  variant of the group
$Z^{2n-2}(X)$, which is suggested by Koll\'ar (cf.
\cite{kollarportug}) given by the following groups:
$$Z^{2n-2}_{rat}(X):=H^{2n-2}(X,\mathbb{Z})/\langle [C],\,C\,\,{\rm rational\,\,curve\,\,in} \,\,X\rangle.$$
Here, by a rational curve, we mean an irreducible curve whose
normalization is rational. These groups are of torsion for $X$
rationally connected, as proved by Koll\'ar (\cite[Theorem 3.13 p 206]{kollarbook}).
It is quite easy to prove that they are birationally invariant.

The proof of Proposition \ref{propdefoinv} gives as well the
following result (already noticed by Koll\'ar  \cite{kollarportug})
:}

\end{rema}
\begin{variant} \label{vardefoinvrat} If $\mathcal{X}\rightarrow B$ is a smooth projective morphism
with rationally connected fibers, the groups $Z^{2n-2}_{rat}(\mathcal{X}_t)$ are
local deformation invariants.
\end{variant}
Let us give one application of Proposition \ref{propdefoinv} (or
rather its proof) and/or its variant \ref{vardefoinvrat}. Let $X$ be
a smooth projective variety of dimension $n+r$, with $n\geq 3$ and
let $\mathcal{E}$ be an ample vector bundle of rank $r$ on $X$. Let
$C_1,\ldots, C_k$ be smooth curves in $X$ whose cohomology classes
generate the group $H^{2n+2r-2}(X,\mathbb{Z})$. For $\sigma\in
H^0(X,\mathcal{E})$, we denote by $X_\sigma$ the zero locus of
$\sigma$. When $\mathcal{E}$ is generated by sections, $X_\sigma$ is
smooth of dimension $n$ for general $\sigma$.
\begin{theo} \label{theocoro}
1) Assume that the sheaves $\mathcal{E}\otimes \mathcal{I}_{C_i}$ are generated by global sections
for $i=1,\ldots,k$. Then if $X_\sigma$ is smooth rationally connected for general $\sigma$,
the group $Z^{2n-2}(X_\sigma)$ vanishes for any $\sigma$ such that $X_\sigma$ is smooth of dimension
$n$.

2) Under the same assumptions as in 1), assume the curves $C_i\subset X$ are rational. Then
if $X_\sigma$ is smooth rationally connected for general $\sigma$,
the group $Z^{2n-2}_{rat}(X_\sigma)$ vanishes for any $\sigma$ such that $X_\sigma$ is smooth of dimension
$n$.
\end{theo}
{\bf Proof.} 1) Let $j_\sigma:X_\sigma\rightarrow X$ be the
inclusion map. Since $n\geq 3$ and $\mathcal{E}$ is ample, by
Sommese's theorem \cite{sommese}, the Gysin map
$j_{\sigma*}:H^{2n-2}(X_\sigma,\mathbb{Z})\rightarrow
H^{2n+2r-2}(X,\mathbb{Z})$ is an isomorphism. It follows that the
group $H^{2n-2}(X_\sigma,\mathbb{Z})$ is a constant group. In order
to show that $Z^{2n-2}(X_\sigma)$ is trivial, it suffices to show
that the classes $(j_{\sigma*})^{-1}([C_i])$ are algebraic on
$X_\sigma$ since they generate $H^{2n-2}(X_\sigma,\mathbb{Z})$.
Since the $X_\sigma$'s are rationally connected, Theorem
\ref{propdefoinv} tells us that it suffices to show that for each
$i$, there exists a $\sigma(i)$ such that $X_{\sigma(i)}$ is smooth
$n$-dimensional and that the class $(j_{\sigma(i)*})^{-1}([C_i])$ is
algebraic on $X_{\sigma(i)}$.

It clearly suffices to exhibit one smooth $X_{\sigma(i)}$ containing $C_i$, which follows
from the following lemma:
\begin{lemm} Let $X$ be a variety of dimension $n+r$ with $n\geq 2$, $C\subset X$ be a smooth curve,
$\mathcal{E}$ be a rank $r$ vector bundle on $X$ such that
$\mathcal{E}\otimes \mathcal{I}_{C}$ is generated by global section. Then for a generic
$\sigma \in H^0(X,\mathcal{E}\otimes \mathcal{I}_{C})$, the zero set $X_\sigma$ is smooth of dimension $n$.
\end{lemm}
{\bf Proof.} The fact that $X_\sigma$ is smooth of dimension $n$ away from $C$ is standard and follows from
the fact that the incidence set
$(\sigma,x)\in \mathbb{P}(H^0(X,\mathcal{E}\otimes \mathcal{I}_{C}))\times (X\setminus C), \,\sigma(x)=0\}$ is smooth
 of dimension
$n+N$, where $N:={\rm dim}\,\mathbb{P}(H^0(X,\mathcal{E}\otimes \mathcal{I}_{C}))$.
It thus suffices to check the   smoothness along $C$ for generic $\sigma$.

This is checked by observing that since $\mathcal{E}\otimes
\mathcal{I}_{C}$ is generated by global sections, its restriction
$\mathcal{E}\otimes N_{C/X}^*$ is also generated by global sections.
This implies that for each point $c\in C$, the condition that
$X_\sigma$ is singular at $c$ defines a codimension $n$ closed
algebraic subset $P_c$ of $P:=\mathbb{P}(H^0(X,\mathcal{E}\otimes
\mathcal{I}_{C}))$, determined by the condition that
$d\sigma_c:N_{C/X,c}\rightarrow \mathcal{E}_c$ is not surjective.
Since ${\rm dim}\,C=1$, the union of the $P_c$'s cannot be equal to
$P$ if $n\geq 2$. \cqfd This concludes the proof of 1) and the proof
of 2) works exactly in the same way. \cqfd

Let us finish this section with the proof of Proposition \ref{spectheo}.

\vspace{0,5cm}

{\bf Proof of Proposition \ref{spectheo}.} Let $p\in{\rm
Spec}\,\mathcal{O}_K$, with residue field $k(p)$. Assume
$\mathcal{X}_p$ is smooth.  For $l$ prime to ${\rm char}\,k(p)$, the (adequately constructed)
specialization map
\begin{eqnarray}
\label{compiso}H^{2n-2}_{et}(X_{\overline{K}},\mathbb{Z}_l(n-1))\rightarrow H^{2n-2}_{et}(\mathcal{X}_{\overline{p}},\mathbb{Z}_l(n-1))\end{eqnarray}
is then an isomorphism (cf. \cite[Chapter VI, {\S}4]{milnebook}).

Observe also that since $X_{\overline{K}}$ is rationally connected, the
rational \'{e}tale  cohomology group $H^{2n-2}_{et}(X_{\overline{K}},\mathbb{Q}_l(n-1))$
is generated over $\mathbb{Q}_l$ by curve classes. Hence the same is true for
$H^{2n-2}_{et}(\mathcal{X}_{\overline{p}},\mathbb{Q}_l(n-1))$. Thus the whole cohomology groups

$$ H^{2n-2}_{et}(X_{\overline{K}},\mathbb{Z}_l(n-1)),\,\, H^{2n-2}_{et}(\mathcal{X}_{\overline{p}},\mathbb{Z}_l(n-1))$$
consist of Tate classes, and (\ref{compiso}) gives an isomorphism
\begin{eqnarray}
\label{compisotate}H^{2n-2}_{et}(X_{\overline{K}},\mathbb{Z}_l(n-1))_{Tate}\rightarrow H^{2n-2}_{et}(\mathcal{X}_{\overline{p}},\mathbb{Z}_l(n-1))_{Tate}.\end{eqnarray}
In order to prove Proposition \ref{spectheo}, it thus suffices to prove the following:

\begin{lemm} 1) For almost every $p\in {\rm Spec}\,\mathcal{O}_K$, the fiber $\mathcal{X}_{\overline{p}}$ is smooth and separably rationally connected.

2) If $\mathcal{X}_{\overline{p}}$ is smooth and separably rationally connected, for any
curve $C_{\overline{p}}\subset \mathcal{X}_{\overline{p}}$, the inverse image $[C_{\overline{p}}]_{\overline{K}}\in H^{2n-2}_{et}(X_{\overline{K}},\mathbb{Z}_l(n-1))$ of the class $[C_{\overline{p}}]\in H^{2n-2}_{et}(\mathcal{X}_{\overline{p}},\mathbb{Z}_l(n-1))$ via the isomorphism
(\ref{compisotate}) is the class of a $1$-cycle  on $X_{\overline{K}}$.
\end{lemm}
{\bf Proof.} 1) When the fiber $\mathcal{X}_p$ is smooth, the
separable rational connectedness of $\mathcal{X}_{\overline{p}}$ is
equivalent to the existence of a smooth rational curve
$C_{\overline{p}}\cong \mathbb{P}^1_{\overline{k(p)}}$ together with
a morphism $\phi:C_{\overline{p}}\rightarrow
\mathcal{X}_{\overline{p}}$ such that the vector bundle
$\phi^*T_{\mathcal{X}_{\overline{p}}}$ on
$\mathbb{P}^1_{{\overline{k(p)}}}$ is a direct sum
$\oplus_i\mathcal{O}_{\mathbb{P}^1_{{\overline{k(p)}}}}(a_i)$ where
all $a_i$ are positive. Equivalently
\begin{eqnarray}\label{vanishing} H^1(\mathbb{P}^1_{{\overline{k(p)}}},\phi^*T_{\mathcal{X}_{\overline{p}}}(-2))=0.
\end{eqnarray}
The smooth projective variety $X_{\overline{K}}$ being rationally
connected in characteristic $0$, it is separably rationally
connected, hence there exists a finite extension $K'$ of $K$, a
curve $C$ and a morphism $\phi:C\rightarrow X$ defined over $K'$,
such that $C\cong \mathbb{P}^1_{K'}$ and
$H^1(\mathbb{P}^1_{K'},\phi^*T_{X_{K'}}(-2))=0$.

We choose a model $$\Phi:\mathcal{C}\cong \mathbb{P}^1_{\mathcal{O}_{K'}}\rightarrow \mathcal{X}'$$
of $C$ and
$\phi$ defined over a Zariski open set of ${\rm Spec}\,\mathcal{O}_{K'}$.
By upper-semi-continuity of cohomology, the vanishing (\ref{vanishing})
remains true after  restriction to almost every closed point $p\in {\rm Spec}\,\mathcal{O}_{K'}$, which proves 1).

2) The proof is identical to the proof of Proposition
\ref{propdefoinv}: we just have to show that the curve
$C_{\overline{p}}\subset\mathcal{X}_{\overline{p}}$ is algebraically
equivalent in $\mathcal{X}_{\overline{p}}$ to a difference
$C''_{\overline{p}}-\sum_i R_{i,{\overline{p}}}$, where each curve
$C''_{\overline{p}}$, resp.  $R_{i,{\overline{p}}}$ (they are in
fact defined over a finite extension $k(p)'$ of $k(p)$), lifts to a
curve $C''$, resp. $R_i$ in $X_{K'}$ for some finite extension $K'$
of $K$.

Assuming the curves $C''_{\overline{p}}$, $R_{i,{\overline{p}}}$ are
smooth, the existence of such a lifting is granted by the condition
$H^1(C''_{\overline{p}},N_{C''_{\overline{p}}/\mathcal{X}_{\overline{p}}})=0$,
resp.
$H^1(R_{i,\overline{p}},N_{R_{i,\overline{p}}/\mathcal{X}_{\overline{p}}})=0$.

Starting from $C\subset \mathcal{X}_{\overline{p}}$ where
$\mathcal{X}_{\overline{p}}$ is separably rationally connected over
$\overline{p}$, we obtain such curves $C''_{\overline{p}}$,
$R_{i,\overline{p}}$ as in the previous proof, applying
\cite[{\S}2.1]{GHS}. \cqfd The proof of Proposition \ref{spectheo} is
finished. \cqfd Again, this proof leads as well to the proof of the
specialization invariance of the $l$-adic analogues
$Z^{2n-2}_{et,rat}(X)_l$ of the groups $Z^{2n-2}_{rat}(X)$
introduced in Remark \ref{remarat}.

\begin{variant} \label{vardefoinv} Let $X$ be a smooth rationally connected variety  defined over a number field
$K$, with ring of integers $\mathcal{O}_K$. Assume given a projective model
$\mathcal{X}$ of $X$ over ${\rm Spec}\,\mathcal{O}_K$. Fix a prime integer $l$. Then
for any  $p\in {\rm Spec}\,\mathcal{O}_K$ such that $\mathcal{X}_{\overline{p}}$ is smooth separably
connected, the group
$Z^{2n-2}_{et,rat}(X)_l$ is isomorphic to the group
$Z^{2n-2}_{et,rat}(X_p)_l$.
\end{variant}
\section{Consequence of a result of Chad Schoen}
In \cite{schoen}, Chad Schoen proves the following theorem:
\begin{theo} \label{theoschoen} Let $X$ be a smooth projective variety of dimension $n$
 defined over a finite field
$k$ of characteristic $p$. Assume that the Tate conjecture holds for
degree $2$ Tate classes on smooth projective surfaces defined over a finite extension of
$k$. Then the \'{e}tale cycle class map:
$$cl : CH^{n-1}(X_{\overline{k}})\otimes \mathbb{Z}_l\rightarrow H^{2n-2}
(X_{\overline{k}},\mathbb{Z}_l(n-1))_{Tate}$$ is surjective, that is
$Z^{2n-2}_{et}(X)_l=0$.
\end{theo}
In other words,  the Tate conjecture \ref{tateconj} for degree $2$
{\it  rational} Tate classes implies that the groups
$Z^{2n-2}_{et}(X)_l$ should be trivial for all smooth projective
varieties defined over finite fields. This is of course very
different from  the situation over $\mathbb{C}$ where the groups
$Z^{2n-2}(X)$ are known to be possibly nonzero.
\begin{rema}{\rm There is a similarity between the proof of Theorem \ref{theoschoen} and the proof of Theorem \ref{theouniruled}. Schoen proves that given an integral Tate class
$\alpha$ on $X$ (defined over a finite field), there exist a smooth
complete intersection surface $S\subset X$ and an integral Tate
class $\beta$ on $S$ such that $j_{S*}\beta=\alpha$ where $j_S$ is
the inclusion of $S$ in $X$. The result then follows from the fact
that if the Tate conjecture holds for degree $2$ rational Tate
classes on $S$, it holds for degree $2$ integral Tate classes on
$S$.

I prove that for $X$ a uniruled or Calabi-Yau, and for $\beta\in
Hdg^4(X,\mathbb{Z})$ there exists surfaces
$S_i\stackrel{j_{S_i}}{\hookrightarrow } X$ (in an adequately chosen
linear system on $X$) and integral Hodge classes $\beta_i\in
Hdg^2(S_i,\mathbb{Z})$ such that $\alpha=\sum_ij_{S_i*}\beta$. The
result then follows from the Lefschetz theorem on $(1,1)$-classes
applied to the $\beta_i$.}

\end{rema}
We refer to \cite{ctz} for some comments on and other applications of Schoen's theorem, and conclude this note with
the proof of the following theorem (cf. Theorem \ref{theoapplischoen} of the introduction).
\begin{theo}\label{theoapplischoenrepete} Assume  Tate's conjecture \ref{tateconj} holds for  degree $2$
 Tate classes on smooth projective surfaces defined over a finite field. Then  the group
$Z^{2n-2}(X)$ is trivial  for  any   smooth rationally connected variety
 $X$ over $\mathbb{C}$.
\end{theo}
{\bf Proof.} We first recall that for a smooth rationally connected variety $X$, the group
$Z^{2n-2}(X)$ is equal to the quotient $H^{2n-2}(X,\mathbb{Z})/H^{2n-2}(X,\mathbb{Z})_{alg}$, due to the fact that
the Hodge structure on $H^{2n-2}(X,\mathbb{Q})$ is trivial. In fact,  we have  more precisely
$$H^{2n-2}(X,\mathbb{Q})=H^{2n-2}(X,\mathbb{Q})_{alg}$$ by hard Lefschetz theorem and the fact that
$$H^{2}(X,\mathbb{Z})=H^{2}(X,\mathbb{Z})_{alg}$$
by the Lefschetz theorem on $(1,1)$-classes.

Next, in order to prove that $Z^{2n-2}(X)$ is trivial, it suffices
to prove that for each $l$, the group $Z^{2n-2}(X)\otimes
\mathbb{Z}_l=H^{2n-2}(X,\mathbb{Z}_l)/({\rm Im}\,cl)\otimes
\mathbb{Z}_l$ is trivial.

We  apply Proposition \ref{propdefoinv} which tells as well that
over $\mathbb{C}$, the group $Z^{2n-2}(X)\otimes \mathbb{Z}_l$ is
locally deformation invariant for families of smooth rationally
connected varieties. Note that our smooth projective rationally
connected variety $X$ is the fiber $X_t$ of a smooth projective
morphism $\phi:\mathcal{X}\rightarrow B$ defined over a number
field, where $\mathcal{X}$ and $B$ are quasiprojective,
geometrically connected and defined over a number field. By local
deformation invariance, the vanishing of $Z^{2n-2}(X)\otimes
\mathbb{Z}_l$ is equivalent to the vanishing of
$Z^{2n-2}(X_{t'})\otimes \mathbb{Z}_l$ for any point $t'\in
B(\mathbb{C})$. Taking for $t'$ a point of $B$ defined over a number
field, $X_{t'}$ is defined over a number field. Hence it  suffices
to prove the vanishing of $Z^{2n-2}(X)\otimes \mathbb{Z}_l$ for $X$
rationally connected defined over a number field $L$.

We have
$$Z^{2n-2}(X)\otimes \mathbb{Z}_l=H^{2n-2}(X,\mathbb{Z}_l)/({\rm Im}\,cl)\otimes \mathbb{Z}_l,$$
and by the Artin comparison theorem (cf. \cite[Chapter III,{\S}3]{milnebook}), this is equal to
$$\frac{H^{2n-2}_{et}(X,\mathbb{Z}_l(n-1))}{({\rm
Im}\,cl)\otimes\mathbb{Z}_l}=Z^{2n-2}_{et}(X)_l$$ since
 $H^{2n-2}_{et}(X,\mathbb{Z}_l(n-1))$ consists of Tate classes. Hence it
suffices to prove that for $X$ rationally connected defined over a
number field and for any $l$, the group $Z^{2n-2}_{et}(X)_l$ is
trivial.

We now apply Proposition \ref{spectheo} to $X$ and its reduction $X_p$
for almost every closed point $p\in {\rm Spec}\,\mathcal{O}_L$. It
follows that the vanishing of $Z^{2n-2}_{et}(X)_l$ is implied by the
vanishing of $Z^{2n-2}_{et}(X_p)_l$. According to Schoen's theorem
\ref{theoschoen}, the last vanishing is implied by the Tate
conjecture for  degree $2$
 Tate classes on smooth projective surfaces.
 \cqfd
\begin{rema}{\rm This argument does not say anything on the groups $Z^{2n-2}_{rat}(X)$, since there is no control on the $1$-cycles representing given degree $2n-2$ Tate classes on
varieties defined over finite fields. Similarly, Theorem \ref{theouniruled} does not say anything on
$Z^4_{rat}(X)$ for $X$ a rationally connected threefold.}
\end{rema}

Institut de math\'{e}matiques de Jussieu

Case 247

4 Place Jussieu

F-75005 Paris, France

\smallskip
 voisin@math.jussieu.fr

\begin{thebibliography}{99}
\bibitem{AH} M. Atiyah, F. Hirzebruch.  Analytic cycles on complex manifolds,
Topology {\bf 1},
25-45 (1962).
\bibitem{B} M. Artin, D. Mumford. Some elementary examples of unirational varieties which are not
rational, Proc. London Math. Soc. (3) 25, 75-95  (1972).
\bibitem{beauville} A. Beauville. Vari\'{e}t\'{e}s rationnelles et unirationnelles, in {\it Algebraic geometry—open problems} (Ravello, 1982), 16-33, Lecture Notes in Math., 997, Springer, Berlin, (1983).
 \bibitem{CG} H. Clemens, P. Griffiths. The intermediate Jacobian of the cubic Threefold, Ann. of Math. 95 (1972), 281-356.
     \bibitem{colliotojanguren} J.-L. Colliot-Th\'{e}l\`{e}ne,  M. Ojanguren. Vari\'{e}t\'{e}s unirationnelles non rationnelles: au-del\`{a} de l'exemple d'Artin et Mumford, Invent. math.  97  (1989),  no. 1, 141--158.
         \bibitem{ctk} J.-L. Colliot-Th\'{e}l\`{e}ne, B. Kahn. Cycles de codimension $2$ et $H^3$ non ramifi\'{e} pour les vari\'{e}t\'{e}s sur les corps finis, preprint 2011.
\bibitem{ctz} J.-L. Colliot-Th\'{e}l\`{e}ne, T. Szamuely. Autour de la conjecture de Tate \`{a} coefficients $\mathbb{Z}_l$ pour les vari\'{e}t\'{e}s sur les corps finis,
 in {\it The Geometry of Algebraic Cycles } (ed.  Akhtar, Brosnan,  Joshua), AMS/Clay Institute Proceedings, 2010, pp. 83-98.
\bibitem{ctv} J.-L. Colliot-Th\'{e}l\`{e}ne, C. Voisin. Cohomologie non ramifi\'{e}e et conjecture
de Hodge enti\`{e}re, preprint 2010 to appear in Duke Math. Journal.
\bibitem{floris} E. Floris. Fundamental divisors on  Fano varieties of index
$n-3$, to appear in Geometriae Dedicata.
\bibitem{GHS} T. Graber, J. Harris, J. Starr. Families of rationally connected varieties.
J. Amer. Math. Soc., 16(1), 57-67  (2003).
 \bibitem{hoeringvoisin} A. H\"{o}ring, C. Voisin. Anticanonical divisors and curve classes on Fano manifolds,  Pure and Applied Mathematics Quarterly
Volume 7, Number 4 (Special Issue: In memory of Eckart Viehweg), 1371-1393 (2011).
\bibitem{kodaira} K. Kodaira. On stability of compact submanifolds of complex manifolds. Amer. J. Math. 85,  79-94 (1963).
\bibitem{kollarbook} J. Koll\'ar. {\it Rational curves on algebraic varieties},  Ergebnisse der Mathematik und ihrer Grenzgebiete. 3. Folge. A Series of Modern Surveys in Mathematics, 32. Springer-Verlag, Berlin, (1996).
 \bibitem{kollar} J. Koll\'ar.  In {\it Trento examples}, Lemma  p. 134
 in {\it Classification of irregular
 varieties}, edited by E. Ballico, F. Catanese, C. Ciliberto,
 Lecture Notes in Math. {\bf  1515}, Springer (1990).
 \bibitem{komimo} J. Koll\'{a}r, Y.  Miyaoka, S. Mori.
   Rationally connected varieties. J. Algebraic Geom. 1  no. 3, 429--448 (1992).
   \bibitem{kollarportug} J. Koll\'ar. Holomorphic and pseudo-holomorphic curves on rationally connected varieties. Port. Math. 67 (2010), no. 2, 155-179.
   \bibitem{milnebook} J. S. Milne. {\it Etale cohomology}, Princeton University Press, 1980.
   \bibitem{milne} J. S. Milne. The Tate conjecture over finite fields (AIM talk),  arXiv:0709.3040.
\bibitem{schoen} C. Schoen. An integral analog of the Tate conjecture for one-dimensional cycles
 on varieties over finite fields. Math. Ann. 311 (1998), no. 3, 493-500.
 \bibitem{sommese} A. Sommese. Submanifolds of abelian varieties, Math. Ann. 233 (1978),
229-256.
\bibitem{soulevoisin} C. Soul\'{e}, C. Voisin.  Torsion cohomology classes and algebraic
 cycles on complex projective manifolds,  Adv. Math.  {\bf  198}  (2005),  no. 1, 107--127.
\bibitem{voisinuniruled} C. Voisin.  On integral Hodge classes on uniruled and Calabi-Yau
threefolds, in {\it Moduli Spaces and Arithmetic Geometry},
Advanced Studies in Pure Mathematics {\bf 45}, 2006, pp. 43-73.
\bibitem{voisinjapjmath} C. Voisin,  Some aspects of the Hodge conjecture,    Japan. J. Math. 2, 261-296 (2007).
    \end{thebibliography}
    \end{document}